\documentstyle[twoside]{article}
\title{A derivation of two quadratic transformations contiguous to that of Gauss via a differential equation approach}
\author{
M. Swathi,\footnote{School of Mathematical and Physical Sciences, Central University of Kerala,
Periye P.O. Dist. Kasaragad 671123, Kerala State, India.
E-Mail: swathimkhd@gmail.com}
\ \ Arjun. K. Rathie\footnote{School of Mathematical and Physical Sciences, Central University of Kerala,
Periye P.O. Dist. Kasaragad 671123, Kerala State, India.
E-Mail: akrathie@cukerala.edu.in} \ \ and 
\ R. B. Paris\footnote{School of Computing, Engineering and Applied Mathematics, University of Abertay Dundee, Dundee DD1 1HG, UK.
E-Mail: r.paris@abertay.ac.uk}\ \footnote{Corresponding author}
 \\}
\begin{document}
\def\f#1#2{\mbox{${\textstyle \frac{#1}{#2}}$}}
\def\dfrac#1#2{\displaystyle{\frac{#1}{#2}}}
\def\boldal{\mbox{\boldmath $\alpha$}}
\newcommand{\bee}{\begin{equation}}
\newcommand{\ee}{\end{equation}}
\newcommand{\lam}{\lambda}
\newcommand{\ka}{\kappa}
\newcommand{\al}{\alpha}
\newcommand{\th}{\theta}
\newcommand{\om}{\omega}
\newcommand{\Om}{\Omega}
\newcommand{\fr}{\frac{1}{2}}
\newcommand{\fs}{\f{1}{2}}
\newcommand{\g}{\Gamma}
\newcommand{\br}{\biggr}
\newcommand{\bl}{\biggl}
\newcommand{\ra}{\rightarrow}
\newcommand{\mbint}{\frac{1}{2\pi i}\int_{c-\infty i}^{c+\infty i}}
\newcommand{\mbcint}{\frac{1}{2\pi i}\int_C}
\newcommand{\mboint}{\frac{1}{2\pi i}\int_{-\infty i}^{\infty i}}
\newcommand{\gtwid}{\raisebox{-.8ex}{\mbox{$\stackrel{\textstyle >}{\sim}$}}}
\newcommand{\ltwid}{\raisebox{-.8ex}{\mbox{$\stackrel{\textstyle <}{\sim}$}}}
\renewcommand{\topfraction}{0.9}
\renewcommand{\bottomfraction}{0.9}
\renewcommand{\textfraction}{0.05}
\newcommand{\mcol}{\multicolumn}
\date{}
\maketitle
\begin{abstract}
The purpose of this note is to provide an alternative proof of two quadratic transformation formulas contiguous to that of Gauss using a differential equation approach.

\vspace{0.4cm}

\vspace{0.3cm}

\end{abstract}

\vspace{0.6cm}

\noindent{\bf 1. \  Introduction}
\setcounter{section}{1}
\setcounter{equation}{0}
\renewcommand{\theequation}{\arabic{section}.\arabic{equation}}
\vspace{0.3cm}

\noindent
The quadratic transformation for the Gauss hypergeometric function ${}_2F_1(a,b;c;x)$ we consider here is the one originally obtained by Gauss in the form (see, for example, \cite[p.~128]{AAR})
\bee\label{e11}
(1+x)^{-2a} {}_2F_1\left[\begin{array}{c} a, b\\2b\end{array}\!;\frac{4x}{(1+x)^2}\right]=
{}_2F_1\left[\begin{array}{c}a, a-b+\fs\\b+\fs\end{array}\!;x^2\right]
\ee
valid when $|x|<1$ and $|4x/(1+x)^2|<1$ and provided $2b$ is neither zero nor a negative integer. 
Bailey \cite{B} re-derived this result by employing the classical Watson summation theorem for the ${}_3F_2$ series. In the standard text of Rainville \cite[p.~63]{R}, the transformation (\ref{e11}) was derived using the differential equation satisfied by ${}_2F_1$.

In 2001, Rathie and Kim \cite{RK} established two transformation formulas contiguous to (\ref{e11}) with the help of a contiguous version of Watson's summation theorem due to Lavoie {\it et al.} These are given in the following theorem.
\newtheorem{theorem}{Theorem}
\begin{theorem}$\!\!\!.$\ \ If $|x|<1$ and $|4x/(1+x)^2|<1$ then

\[(1+x)^{-2a} {}_2F_1\left[\begin{array}{c} a, b\\2b+1\end{array}\!;\frac{4x}{(1+x)^2}\right]=
{}_2F_1\left[\begin{array}{c}a, a-b+\fs\\b+\fs\end{array}\!;x^2\right]\hspace{3cm}\]
\bee\label{e12}
\hspace{4cm}-\frac{2ax}{2b+1}\,\,{}_2F_1\left[\begin{array}{c}a+1, a-b+\fs\\b+\f{3}{2}\end{array}\!;x^2\right]
\ee
and
\[(1+x)^{-2a} {}_2F_1\left[\begin{array}{c} a, b\\2b-1\end{array}\!;\frac{4x}{(1+x)^2}\right]=
{}_2F_1\left[\begin{array}{c}a, a-b+\f{3}{2}\\b-\fs\end{array}\!;x^2\right]\hspace{3cm}\]
\bee\label{e13}
\hspace{4cm}+\frac{2ax}{2b-1}\,\,{}_2F_1\left[\begin{array}{c}a+1, a-b+\f{3}{2}\\b+\fs\end{array}\!;x^2\right]
\ee
provided $2b\pm 1$ is neither zero nor a negative integer, respectively.
\end{theorem}
Here we give an alternative demonstration of the quadratic transformations (\ref{e12}) and (\ref{e13}) by adopting the differential equation approach employed by Rainville. It is worth remarking that these transformations cannot be derived completely by the hypergeometric differential equation, but that a related second-order differential equation has to be solved by the standard Frobenius method.

Before we give our alternative derivation of (\ref{e12}) and (\ref{e13}) in Section 3, we first present an outline of the arguments employed by Rainville \cite[p.~63]{R} to establish the Gauss transformation (\ref{e11}).

\vspace{0.6cm}

\noindent{\bf 2. \ Derivation of (\ref{e11}) by Rainville's method}
\setcounter{section}{2}
\setcounter{equation}{0}
\renewcommand{\theequation}{\arabic{section}.\arabic{equation}}
\vspace{0.3cm}

\noindent 
The hypergeometric function ${}_2F_1(a,b;c;z)$ satisfies the differential equation \cite[p.~75]{AAR}, \cite[Eq.~(15.10.1)]{DLMF}
\bee\label{e21}
z(1-z) \frac{d^2w}{dz^2}+\{c-(a+b+1)z\} \frac{dw}{dz}-abw=0.
\ee
If we put $c=2b$ and make the change of variable $z=4x/(1+x)^2$, then equation (\ref{e21}) becomes
\[x(1-x)(1+x)^2 \frac{d^2w}{dx^2}+2(1+x)\{b-2ax+(b-1)x^2\} \frac{dw}{dx}-4(1-x)abw=0.\]
If we now put $w=(1+x)^{2a}y$, then after some simplification we find
\bee\label{e22}
x(1-x^2) \frac{d^2y}{dx^2}+2\{b-(2a-b+1)x^2\} \frac{dy}{dx}-2ax(1+2a-2b)y=0,
\ee
of which one solution is
\bee\label{e23}
y=(1+x)^{-2a}\,{}_2F_1\left[\begin{array}{c} a, b\\2b\end{array}\!;\frac{4x}{(1+x)^2}\right].
\ee

The differential equation (\ref{e22}) is invariant under the change of variable from $x$ to $-x$. Hence, if we introduce the new independent variable $v=x^2$ the equation describing $y$ becomes
\bee\label{e24}
v(1-v) \frac{d^2y}{dv^2}+\{b+\fs-(2a-b+\f{3}{2})v\} \frac{dy}{dv}-a(a-b+\fs)y=0.
\ee
We observe that (\ref{e22}) is of the same form as the hypergeometric differential equation (\ref{e21}), which therefore has in $|v|<1$ the two solutions \cite[Eq.~(15.10.2)]{DLMF}
\bee\label{e25}
 {}_2F_1\left[\begin{array}{c} a, a-b+\fs\\b+\fs\end{array}\!;v\right]\quad\mbox{and}\quad
v^{\frac{1}{2}-b}{}_2F_1\left[\begin{array}{c} a-b+\fs, a-2b+1\\\f{3}{2}-b\hspace{0.2cm}\end{array}\!;v\right].
\ee
We observe that the differential equation (\ref{e22}) has the solution (\ref{e23}) valid in $|4x/(1+x)^2|<1$,
provided $2b$ is neither zero nor a negative integer. At the same time, equation (\ref{e22}) has the solutions (\ref{e25}) with $v=x^2$ valid in $|x|<1$. Therefore, subject to these conditions, there exist constants $A$ and $B$ such that
\[(1+x)^{-2a} {}_2F_1\left[\begin{array}{c} a, b\\2b\end{array}\!;\frac{4x}{(1+x)^2}\right]=
A\,{}_2F_1\left[\begin{array}{c}a, a-b+\fs\\b+\fs\end{array}\!;x^2\right]\hspace{2cm}\]
\[\hspace{5.2cm}+Bx^{1-2b} {}_2F_1\left[\begin{array}{c}a-b+\fs, a-2b+1\\\f{3}{2}-b\hspace{0.2cm}\end{array}\!;x^2\right].\]
The left-hand side and the first member on the right-hand side of the above expression are both analytic at $x=0$, but the remaining term is not due to the presence of the factor $x^{1-2b}$. Hence $B=0$ and by considering the terms at $x=0$ it is easily seen that $A=1$. This leads to the required quadratic transformation given in (\ref{e11}).

\vspace{0.6cm}

\noindent{\bf 3. \ An alternative derivation of Theorem 1}
\setcounter{section}{3}
\setcounter{equation}{0}
\renewcommand{\theequation}{\arabic{section}.\arabic{equation}}
\vspace{0.3cm}

\noindent 
We first establish the quadratic transformation (\ref{e12}). With $c=2b+1$ in (\ref{e21}) and the change of variable $z=4x/(1+x)^2$ we obtain
\[x(1-x)(1+x)^2 \frac{d^2w}{dx^2}+(1+x)\{2b+1-4ax+2x+(2b-1)x^2\} \frac{dw}{dx}-4ab(1-x)w=0,\]
which has a solution $w={}_2F_1(a,b;2b+1;4x/(1+x)^2)$. With the further change of dependent variable
$w=(1+x)^{2a}y$, we find after some simplification
\bee\label{e31}
x(1-x^2) \frac{d^2y}{dx^2}+\{2b+1+2x-(4a-2b+1)x^2\} \frac{dy}{dx}+2a\{1+2(b-a)x\} y=0.
\ee
A solution of (\ref{e31}) is consequently
\[y=(1+x)^{-2a} \,{}_2F_1\left[\begin{array}{c} a, b\\2b+1\end{array}\!;\frac{4x}{(1+x)^2}\right].\]

The differential equation (\ref{e31}) is not invariant under the change of variable $x$ to $-x$, and so we cannot reduce it to the hypergeometric equation (\ref{e21}). Inspection of (\ref{e31}) shows that the point $x=0$ is a regular singular point. Accordingly, we seek two linearly independent solutions of (\ref{e31}) by the Frobenius method and let
\bee\label{e32}
y=x^\lambda \sum_{n=0}^\infty c_nx^n\qquad (c_0\neq 0),
\ee
where $\lambda$ is the indicial exponent. Substitution of this form for $y$ in (\ref{e31}) then leads after a little simplification to
\[\sum_{n=0}^\infty c_nx^{n-1}(n+\lambda)(n+\lambda+2b)=
\sum_{n=0}^\infty c_nx^{n+1}(n+\lambda+2a)(n+\lambda+2a-2b)\]
\[\hspace{5cm}-2\sum_{n=0}^\infty c_nx^n (n+\lambda+a).\]
The coefficients of $x^{-1}$ must vanish to yield the indicial equation
\[ \lambda(\lambda+2b)=0,\]
so that $\lambda=0$ and $\lambda=-2b$. Equating the coefficients of $x^n$ for non-negative integer $n$, we obtain
\begin{eqnarray}
c_1&=&\frac{-2(\lambda+a)}{(1+\lambda)(1+\lambda+2b)}\,c_0,\nonumber\\
\label{e33}\\
c_{n}&=&\frac{\{n+\lambda+2(a-1)\}\{n+\lambda+2(a-b-1)\}c_{n-2}-2(n+\lambda+a-1)c_{n-1}}{(n+\lambda)(n+\lambda+2b)}
\quad(n\geq 2).\nonumber
\end{eqnarray}

With the choice $\lambda=0$, we have
\begin{eqnarray*}
c_1&=&\frac{-2a}{(2b+1)}\,c_0,\\
c_{n}&=&\frac{\{n+2(a-1)\}\{n+2(a-b-1)\}c_{n-2}-2(n+a-1)c_{n-1}}{n(n+2b)}
\qquad(n\geq 2).
\end{eqnarray*}
Solution of this three-term recurrence with the help of {\it Mathematica} generates the values given by
\[c_{2n}=\frac{(a)_n (a-b+\fs)_n}{n!\, (b+\fs)_n}\,c_0,\qquad c_{2n+1}=\frac{(a+1)_n (a-b+\fs)_n}{n!\, (b+\f{3}{2})_n}\,c_1,\]
the general values being established by induction. Substitution in (\ref{e32}) then yields one solution of (\ref{e31}) given by
\[y_1=c_0\left\{{}_2F_1\left[\begin{array}{c}a, a-b+\fs\\b+\fs\end{array}\!;x^2\right]-\frac{2ax}{2b+1}\,{}_2F_1\left[\begin{array}{c}a+1, a-b+\fs\\b+\f{3}{2}\end{array}\!;x^2\right]\right\}\]
when $|x|<1$.

A second solution is obtained by taking $\lambda=-2b$ in (\ref{e33}) to yield
\begin{eqnarray*}
c_1&=&\frac{-2(a-2b)}{(1-2b)}\,c_0,\\
c_{n}&=&\frac{\{n+2(a-2b-1)\}\{n+2(a-b-1)\}c_{n-2}-2(n+a-2b-1)c_{n-1}}{n(n-2b)}
\qquad(n\geq 2).
\end{eqnarray*}
This generates the values
\[c_{2n}=\frac{(a-2b)_n (a-b+\fs)_n}{n!\, (\fs-b)_n}\,c_0,\qquad c_{2n+1}=\frac{(a-2b+1)_n (a-b+\fs)_n}{n!\, (\f{3}{2}-b)_n}\,c_1.\]
A second solution of (\ref{e31}) is therefore given by
\[y_2=c_0x^{-2b}\left\{{}_2F_1\left[\begin{array}{c}a-2b, a-b+\fs\\\fs-b\end{array}\!;x^2\right]-\frac{2(a-2b)x}{1-2b}\,{}_2F_1\left[\begin{array}{c}a-2b+1, a-b+\fs\\\f{3}{2}-b\end{array}\!;x^2\right]\right\}\]
when $|x|<1$.

It then follows, when $|x|<1$ and $|4x/(1+x)^2|<1$ and provided $2b+1$ is neither zero nor a negative integer,
that there exist constants $A$ and $B$ such that
\bee\label{e34}
(1+x)^{-2a} {}_2F_1\left[\begin{array}{c} a, b\\2b+1\end{array}\!;\frac{4x}{(1+x)^2}\right]=Ay_1+By_2.
\ee
Now the left-hand side of (\ref{e34}) and the solution $y_1$ are both analytic at $x=0$, whereas the solution $y_2$ is not
analytic at $x=0$ due to the presence of the factor $x^{-2b}$. Hence $B=0$ and, by putting $x=0$ in (\ref{e34}), it is easily seen that $A=1$. This then yields the result stated in (\ref{e12}).

A similar procedure can be employed to establish the quadratic transformation in (\ref{e13}). Putting $c=2b-1$ in (\ref{e21}) and carrying out the same sequence of transformations, we obtain the differential equation satisfied by 
\bee\label{e35}
y=(1+x)^{-2a} \,{}2F_1\left[\begin{array}{c}a, b\\2b-1\end{array}\!;\frac{4}{(1+x)^2}\right]
\ee
in the form
\bee\label{e36}
x(1-x^2) \frac{d^2y}{dx^2}+\{2b-1-2x-(4a-2b+3)x^2\} \frac{dy}{dx}-2a\{1+2(a-b+1)x\}y=0.
\ee
Substitution of (\ref{e32}) then leads to the three-term recurrence for the coefficients $c_n$
\begin{eqnarray*}
c_1&=&\frac{-2(\lambda+a)}{(1+\lambda)(\lambda+2b-1)}\,c_0,\\
c_{n}&=&\frac{\{n+\lambda+2(a-1)\}\{n+\lambda+2(a-b)\}c_{n-2}+2(n+\lambda+a-1)c_{n-1}}{(n+\lambda)(n+\lambda+2b-2)}
\quad(n\geq 2)
\end{eqnarray*}
subject to the indicial equation $\lambda(\lambda+2b-2)=0$. The choice of indicial exponent $\lambda=0$ yields with the help of {\it Mathematica} the values of the coefficients given by
\[c_{2n}=\frac{(a)_n (a-b+\f{3}{2})_n}{n!\, (b-\fs)_n}\,c_0,\qquad c_{2n+1}=\frac{(a+1)_n (a-b+\f{3}{2})_n}{n!\, (b+\fs)_n}\,c_1,\]
with $c_1=2a/(2b-1)$, and the choice $\lambda=2-2b$ yields
\[c_{2n}=\frac{(a-b+2)_n(a-b+\f{3}{2})_n }{n!\, (\f{3}{2}-b)_n}\,c_0,\qquad c_{2n+1}= \frac{(a-2b+3)_n(a-b+\f{3}{2})_n }{n!\, (\f{5}{2}-b)_n}\,c_1,\]
with $c_1=2(a-2b+2)/(3-2b)$.

Consequently two solutions of the differential equation (\ref{e36}) are
\[
y_1=c_0\left\{{}_2F_1\left[\begin{array}{c}a, a-b+\f{3}{2}\\b-\fs\end{array}\!;x^2\right]+\frac{2ax}{2b-1}\,{}_2F_1\left[\begin{array}{c}a+1, a-b+\f{3}{2}\\b+\fs\end{array}\!;x^2\right]\right\}\]
and
\[y_2=c_0x^{2-2b}\left\{{}_2F_1\left[\begin{array}{c}a-b+2, a-b+\f{3}{2}\\\f{3}{2}-b\end{array}\!;x^2\right]\right.\hspace{4cm}\]
\[\left.\hspace{4cm}+\frac{2(a-2b+2)x}{3-2b}\,{}_2F_1\left[\begin{array}{c}a-2b+3, a-b+\f{3}{2}\\\f{5}{2}-b\end{array}\!;x^2\right]\right\}\]
when $|x|<1$. It then follows, when $|x|<1$, $|4x/(1+x)^2|<1$ and provided $2b-1$ is neither zero nor a negative integer, that there exist constants $A$ and $B$ such that the function in (\ref{e35}) can be expressed as $Ay_1+By_2$. For the same reasons as in the previous case we find $A=1$ and $B=0$, thereby establishing (\ref{e13}).

\end{document}